\documentclass{amsart}

\usepackage{amssymb}
\usepackage{amscd}
\usepackage[all]{xy}

\DeclareMathOperator{\Hom}{{\rm Hom}}

\DeclareMathOperator{\Fr}{\rm Fr}
\DeclareMathOperator{\Gal}{\rm Gal}

\DeclareMathOperator{\cor}{\rm Cor}
\DeclareMathOperator{\sgn}{\rm sgn}

\newcommand{\bbz}{\mathbf z}
\newcommand{\ZZ}{\mathbb Z}

\newcommand{\QQ}{{\mathbb Q}}

\newcommand{\T}{{\mathcal T}}

\newcommand{\AAA}{{\mathcal A}}
\newcommand{\BBB}{{\mathcal B}}
\newcommand{\DDD}{{\mathcal D}}
\newcommand{\LLL}{{\mathcal L}}
\newcommand{\NNN}{{\mathcal U}}
\newcommand{\co}{{\mathcal O}}

\newcommand{\K}{\mathbf K}

\newcommand{\Kc}{\mathbf K^{\bullet,\bullet,\bullet}}
\newcommand{\Q}{\mathbf Q}

\newcommand{\yf}{Y_{\text{fin}}}
\newcommand{\zf}{Z_{\text{fin}}}

\newtheorem{theo}{Theorem}[section]
\newtheorem{Lemma}[theo]{Lemma}
\newtheorem{Proposition}[theo]{Proposition}

\theoremstyle{remark}
\newtheorem{rem}[theo]{Remark}

\title{On the universal norm distribution}
\subjclass{Primary 11G99, 11R23; Secondary 11R18}
\date{October 14, 2002}
\author{Yi Ouyang}
\address{
Department of Mathematics, University of Toronto, 100 St. George St.,
Toronto, ON M5S 3G3, Canada}
\email{youyang@math.toronto.edu}

\begin{document}
\begin{abstract}
We introduce and study the universal norm distribution in this paper,
which generalizes the concepts of universal ordinary distribution
and the universal Euler system. We study the Anderson type resolution
of the universal norm distribution and then use this resolution
 to study the group cohomology of the universal norm distribution.
\end{abstract}

\maketitle
\section{Introduction}
Let $r$ be a positive integer, the \emph{universal ordinary distribution} 
of rank $1$ and level $r$ is well known to be the free abelian group
\[ U_r=\frac{\langle [a]:a\in \frac{1}{r}\ZZ/\ZZ\rangle}
{\langle [a]-\sum_{pb=a}[b]: p\mid r, a\in\frac{p}{r}\ZZ/\ZZ \rangle}. \]
With a natural $G_r=\Gal(\QQ(\mu_r)/\QQ)$ action on $U_r$, 
$U_r$ becomes a $G_r$-module and 
plays a very important role in the study of cyclotomic fields,
see for example Lang~\cite{Lang} or Washington~\cite{Washington}
for more information. 
In particular, the sign cohomology of $U_r$ gives key information 
about the indices of cyclotomic units and Stickelberger ideals
as illustrated by Sinnott's original paper~\cite{Sinnott} and many
following papers on this subject by different authors. 
The $G_r$-cohomology is found to be related to the cyclotomic 
Euler system, as shown by Anderson-Ouyang~\cite{AO} about
the Kolyvagin recursion in the universal ordinary distribution.

In the book~\cite{Rubin2}, Rubin defined a generalization of
the universal ordinary distribution, which he called the 
\emph{universal Euler system}. It then was used to prove the 
Kolyvagin recursions satisfied by the Euler systems.  
However, there are other 
universal objects satisfying similar distribution relations. In the
paper~\cite{Ouyang3}, we proposed a generalization of the universal
ordinary distribution, for which we called the \emph{universal norm
distribution}. We used it successfully to study Sinnott's index formula.

We further generalize the idea of the universal norm distribution 
in this paper, which  treats the universal Euler systems 
as special cases. We study in detail the structure of the universal
norm distribution in this paper. We also study in detail its group 
cohomology. In short, this paper generalizes
the results of Ouyang~\cite{Ouyang2} and the appendix of it by Anderson.
The goal is to set up necessary tools to the study the universal
Kolyvagin recursion for the universal norm distribution(thus includes
the universal Euler system case), which is a question raised in
Anderson-Ouyang~\cite{AO} and will be answered  in a
subsequent paper~\cite{Ouyang4}. However, our study here is more than
applications to the universal Kolyvagin recursion. The pure
homological setup here should offer us more freedom to the study
of other arithmetic aspects of the universal norm distribution. 
Certainly we expect more studies in this direction.

The structure of this paper is as follows.
We first introduce the definition of the universal norm distribution
$\NNN_z$ in \S~\ref{section:not} and give some examples
in \S~\ref{section:exp}. Basic properties of $\NNN_z$ is studied
in \S~\ref{section:basic}.
A general phenomenon of  every universal norm distribution $\NNN_z$
is  Anderson's resolution $\LLL_z$ attached 
to it. We construct $\LLL_z$ in \S~\ref{section:and} and prove
it is indeed a resolution of $\NNN_z$ in
Theorem~\ref{Proposition:Anderson}, a generalization of the results
by Anderson in the appendix of \cite{Ouyang2}.  Because of the existence
of Anderson's resolution $\LLL_z$, 
we can thus  apply the double complex and spectral sequences method 
to study the group cohomology of the universal 
norm distribution $\NNN_z$. This is accomplished in \S~\ref{sec:ggt},
in particular, in Theorem~\ref{theo:b} and Theorem~\ref{theo:a}.
For the universal ordinary distribution
case, the two Theorems recover and generalize Theorem A
in Ouyang~\cite{Ouyang2}.

The author got very first idea of this paper during his pleasant 
visit in IHES in Spring 2001. Part of the results here was reported in 
the number theory
seminar in Penn State University in November 2001 and then in McMaster
University in February 2002, and in the summer meeting of CMS at Laval
University in June 2002. The author sincerely thanks the above 
organizations, Professors
Robert Vaughn and Winnie Li at PSU, Professor Manfred Kolster
at McMaster and Professors Kumar Murty and Paramath Sastry at Toronto 
for inviting me to give these talks. Last but not least, thanks always go to 
Professor Greg W. Anderson for his ideas and his influence.

\section{Notations and Definitions} \label{section:not}
\subsection{Basic Notations} \label{subsection:Notations}
Let $X$ be a totally ordered set. Denote by $x$, $x_i$ the elements 
in $X$.  

Let $Y$ be the set of all squarefree formal products of $x\in X$, i.e.,
the element $y\in Y$ has the form $x_1\cdots x_n\cdots $ for 
$x_i\neq x_j\in X$. In particular, let ${\mathbf 1}\in Y$ denote the 
element of which no $x\in X$ appears in the formal product.
One can identify $Y$ with the collection of all subsets of $X$, thus
$\mathbf 1$ is corresponding to the empty set. For every $y\in Y$, the 
\emph{degree} $\deg y$ of $y$ is define to be the number of elements 
$x\in X$ dividing $y$. Denote by $y$, $y_i$ the 
elements in $Y$. If without further statement, we'll assume that $y$ is 
finite, i.e., $\deg y<\infty$. Denote by $\yf$ the set of all finite
$y\in Y$. 

Let $Z$ be the set of all formal product of $x\in X$, 
i.e., the element
$z\in Z$ has the form $x^{i_1}_1\cdots x^{i_n}_n\cdots$ with
$i_j\in \ZZ_{\geq 0}$. For every 
$z=x^{i_1}_1\cdots x^{i_n}_n\cdots$, define the \emph{degree} of $z$
to be $\deg z=\sum_{j=1}^n i_j$.
Denote by $z$, $z'$, $w$ the elements in $Z$ and in particular by 
$\bbz$, $\bbz'$ the infinite elements(elements with infinite degree) 
in $Z$.  The subset of all finite elements in $Z$ will be denoted by $\zf$.

Apparently we have $X\subseteq Y\subseteq Z$. One can always keep in 
mind the example that $X$ is the set of prime numbers, $\yf$ is the set of all
squarefree positive integers and $\zf$ is the set of positive integers.
We can thus imitate all the terminologies from traditional sense, for 
example, prime factors, factors, the greatest common divisors and etc.

For every $z\in Z$ and $x\in X$, the \emph{valuation} of $z$ at $x$ is
the highest power of $x$ dividing $z$ and is denoted by $v_x(z)$.
For every $z\in Z$, there exists a unique $\bar z\in Y$($\bar z$ could be
infinite) such that 
if $x\mid z$ then $x\mid \bar z$. We call $\bar z$ the \emph{support} 
of $z$. For every $z\in Z$, if a factor $z'\mid z$ satisfies 
$\gcd (z', z/z')=1$,  $z'$ is called a \emph{stalk} of $z$ and 
is denoted by $z'\mid_s z$. 
Note that the set of stalks of $z$ has a one-to-one correspondence with the 
set of factors(and also stalks) of $\bar z$.
Fix $z$, for each $y\mid \bar z$, let $z(y)$ be the stalk of $z$ whose 
support is $y$. In particular, $z(x)$ is just $x^{v_x(z)}$.

For each pair $x\in X$ and $y\in Y$, we define the function
$\omega: \ X\times Y\rightarrow \{1,0, -1\}$ by
\[ (x,y)\longmapsto \begin{cases} (-1)^{\#\{x': x'<x\}},\ 
&\text{if $x\mid y$};\\ 0,\ 
&\text{if $x\nmid y$.} \end{cases} \]

Let $G$ be a profinite group. Let $A$ be a point set with discrete topology
such that $G$ acts continuously. Suppose there is a surjection
$A\rightarrow \zf$ which induces a bijection between the orbits of $A$
and elements $z\in\zf$. Let $B_z$ be the corresponding orbit of $z$. Let 
$H_z$ be the stabilizer of any $b\in B_z$. 
We assume $\{H_z: z\in \zf\}$ satisfies the following axioms:
\begin{itemize}
\item For every $z\in \zf$, the commutator $[G,G]\leq H_z$;
\item  For every $z'\mid z\in \zf$, $H_{z}\leq H_{z'}$;
\item For $z$ and $z'$ in $\zf$ and relatively prime, 
$H_{zz'}=H_z \cap H_{z'}$ and $G=H_z H_{z'}$.       
\end{itemize} 
By the first axiom, then $H_z$ is a normal
open subgroup of $G$ and the quotient group
$G_z=G/H_z$ is finite abelian.
By the second axiom, for every $z'\mid z\in \zf$, $G_{z'}$ is a 
quotient group of $G_{z}$;
by the last axiom, one see that for every $z'\mid_s z$, the 
quotient map $G_z\rightarrow G_{z'}$ is 
canonically split as $G_z=G_{z'}\times G_{z/z'}$, we thus have the 
following canonical decomposition
\[ G_z=\prod_{x\mid z} G_{z(x)}. \]
Let $N_z$ be the sum of all 
elements $g\in G_z$ in the group ring $\ZZ[G_z]$. 
For $z$ finite and $z'\mid z$, 
Let $g_{z'}$ denote the image of $g\in G_z$ in $G_{z'}$. Let
$N^{z}_{z'}$ be the corresponding inflation map from $\ZZ[G_{z'}]$ to 
$\ZZ[G_z]$.  For every infinite
$\bbz\in Z$, let $G_{\bbz}$ be the inverse limit of $G_z$ over all finite
$z\mid_s \bbz$. Then  $G_{\bbz}$ is actually the direct product of 
$G_{\bbz(x)}$ for every $x\mid \bbz$.  

Write $B_z=\{[gz]: g\in G_z\}$, then 
\[ A=\bigcup_{z\in \zf} B_z=\{[gz]: g\in G_z, z\in \zf\}, \]
and $G_{x^n}$ acts trivially in $B_z$ if $x\nmid z$. Thus $A$ and
$\{G_z: z\in \zf\}$ are uniquely determined by each other.  
Let $A_{z}=\bigcup_{z'\mid_s z, z'\in \zf} B_{z'}$ for every $z\in Z$.

For each pair $x\in X$ and $z\in Z$, the Frobenius element
$\Fr_{x}$ is a given element in $G$ whose restriction to $G_{x^n}$ is the 
identity for every $n\in \mathbb N$.

Let $\co$ be an integral domain and let $\Phi$ be its fractional field. 
Let $\T$ be a fixed $\co$-algebra which is torsion free and 
finitely generated as an $\co$-module.
We suppose that $\T$ is a trivial $G$-module. For each $x\in X$, a polynomial 
\[ p(x; t)\in \T[\, t\,] \] 
is chosen corresponding to $x$.

\subsection{Definition of the universal norm distribution} 
Let $\AAA$ be the free $\T$-module generated by $A$, along with the
$G$-action, $\AAA$ becomes a torsion free $\T[G]$-module. Let
$\BBB_z$ be the $\T[G]$-submodule of $A$ generated by $B_z$ as $\T$-module
for $z\in \zf$.
Then $\BBB_z$ is nothing but a free rank $1$ $\T[G_z]$-module with generator
$[z]$. Let $\AAA_z$ be the $\T[G]$-submodule generated by $A_z$ 
as $\T$-module for every
$z\in Z$. Thus $\AAA_z$ has a natural $\T[G_{z'}]$-module structure
for every $z\mid_s z'$.

Let $\lambda_{z(x)}$
be the $\T[G_z]$-homomorphism of $\AAA_z$ given by 
\[ \lambda_{z(x)}: [z']\longmapsto \begin{cases} 
p(x; \Fr^{-1}_{x})[z']-N_{z(x)}[z(x)z'],\
&\text{if}\ x\nmid z',\\ 0,\ &\text{if}\ x\mid z'.\end{cases} \]
Let $\DDD_z$ be the submodule of $\AAA_z$ generated by the images
of $\lambda_{z(x)}(\AAA_{z/z(x)})$ for all $x\mid z$. Elements in $\DDD_z$ are 
called \emph{distribution relations} in $\AAA_z$.
The \emph{universal norm distribution}
$\NNN_z$ according to the above assumptions is defined to be the
quotient $\T[G_z]$-module $\AAA_z/\DDD_z$, i.e., $\AAA_z$ modulo all
distribution relations.

Note that for every $z\in Z$,
\[ \AAA_z=\bigcup_{\substack{z'\ \text{finite}\\ z'\mid_s z}}\AAA_{z'}. \]
For any $z'\mid_s z$, the apparent inclusion of $\AAA_{z'}$ to 
$\AAA_z$ induces an injection map from 
$\NNN_{z'}$ to $\NNN_z$. In Proposition~\ref{Proposition:basis2}(2), we'll
see this injection actually is a splitting $G_z$-monomorphism.

\section{Examples} \label{section:exp}
We give a few examples about the universal norm distribution here.

\subsection{The trivial case} The first case of the universal norm 
distribution is that $p(x;t)=1$ for every $x\in X$. In this case, one 
easily see that $\NNN_z$ is generated by the images of $B_z$. Actually,
$\NNN_z$ is nothing but isomorphic to the $\T$-module $\BBB_z=\T[G_z]$
(see the remark after Proposition~\ref{Proposition:basis2}). 
We call this type of universal norm distribution
the \emph{trivial universal norm distribution}.

\subsection{The universal ordinary distribution} \label{subsection:ord}
Recall that an \emph{ordinary
distribution} of level $r$ for a positive integer $r$ is a function 
$f$ from  $\frac{1}{r}\ZZ/\ZZ$ to an abelian group $Ab$ satisfying
\[ f(pa)=\sum_{i=0}^p f(a+\frac{i}{p}),\ \forall\ \text{primes}\ p\mid r. \]
In the category of ordinary distributions, there exists a universal object,
i.e., an abelian group $U_r$ and a distribution relation 
$u: \frac{1}{r}\ZZ/\ZZ \rightarrow U_r$ such that for every $f$, there is a
unique homomorphism $f^u: U_r\rightarrow Ab$, such that $f=f^u\circ u$. 
Usually one can write $U_r$ as
\[ \frac{\langle [a]: a\in \frac{1}{r}\ZZ/\ZZ\rangle}
{\langle [pa]-\sum_{i=0}^p [a+\frac{i}{p}]: \forall p\mid r\rangle} \]
and the map $u$ sends $a$ to $[a]$.

The universal ordinary distribution $U_r$ is actually a universal norm
distribution according to our language. 
Let $X$ be the set of all prime numbers. Then $\yf$ is the set of all
squarefree positive integers and 
$\zf$ is just the set of positive integers. Let $G=G_{\QQ}$. 
Let $G_r=\Gal(\QQ(\zeta_{r})/\QQ)$. The Frobenius element
$\Fr_{p}$ is defined by the usual way. Let  $\co=\ZZ=\T$ 
and thus $\Phi=\QQ$. Let the polynomial $p(p; t)=1-t$ for all $p\in X$. 
The corresponding universal norm distribution $\NNN_{r}$ is shown 
to be isomorphic to the universal
ordinary distribution $U_r$(see Ouyang~\cite{Ouyang3}) by sending
$[r]\in \NNN_r$ to $[\frac{1}{r}]\in U_r$.

\subsection{The universal ordinary predistribution}
Keep $X$, $Y$, $Z$, $G$, $\co$ and $\T$ the same as in \S~\ref{subsection:ord}.
Now let $p(p;t)=-t$ for $p\neq 2$ and let $p(2;t)=-t$, we call the
resulting universal norm distribution the  \emph{universal ordinary
predistribution}. 
\begin{Proposition} The universal ordinary predistribution
is isomorphic to the integer ring of the cyclotomic
number field $\QQ(\mu_r)$ for each $r$.
\end{Proposition}
\begin{proof} Define 
$\mathbf e_r: \AAA_r\rightarrow \co_{\QQ(\mu_r)}$
by
\[ [\sigma r']\longmapsto \exp(\frac{2\pi i}{r'})^{\sigma},  \]
then immediately one has
\begin{enumerate}
\item $\DDD_r\subseteq \ker \mathbf e_r$,
\item $e_r$ is surjective.
\end{enumerate}
By Proposition~\ref{Proposition:basis} which we'll prove later, we know
that $\NNN_r$ has $\ZZ$-rank $\phi(r)$, the same as $\co_{\QQ(\nu_r)}$,
thus $\mathbf e_r$ is an isomorphism.
\end{proof}

\subsection{The universal Euler system} 
Let $K$ be  a fixed number field. Let $p$ be a rational prime number.
Let  $\Phi$ be a finite extension of $\QQ_p$ and let $\co$ be the
ring of integer of $\Phi$. Let $T$ be a $p$-adic representation of
$G_K$ with coefficients in $\co$. Assume that $T$ is unramified outside a 
finite set of primes of $K$. 

Fix an ideal $\mathfrak N$ of $K$ divisible by $p$ and by all primes where 
$T$ is ramified. Let $X$ be the set of all primes $x$ of $K$ which is
 prime to $\mathfrak N$
and $K(x)\neq  K({\mathbf 1})$, where $K(x)$ is the ray class field
of $K$ modulo $x$. Then $Y$ and $Z$ are defined following $X$. For
every $y=x_1\cdots x_n\in Y$, let $K(y)$ be the 
composite
\[ K(y)=K(x_1)\cdots K(x_n). \]
Fix a $\ZZ^d_p$-extension $K_{\infty}/K$ which no finite prime
splits completely. We write $K\subset_f F\subset K_{\infty}$ to indicate
$F/K$ a finite subextension of $K_{\infty}/K$. For 
$K\subset_f F\subset K_{\infty}$, we let
$F(y)=FK(y)$. Let $G_y=\Gal(F(y)/F(\mathbf 1))\cong \Gal(K(y)/K(\mathbf 1))$.
We see that for any $y'\mid y$, $G_{y}=G_{y'}\times G_{y/y'}$. 
Let $G=G_{K(\mathbf 1)}$. 

Let $\Fr_{x}$ denote a Frobenius of $x$ in $G_K$, and let
\[ p(x;t)=\det(1-\Fr_{x}^{-1} t| T^{\ast}) 
\in \co[t]. \]

Let $\T=\co[\Gal(F(\mathbf 1)/K)]$. With the above $X$, $Y$, $\co$, $\Phi$
and $p(x;t)$, the corresponding universal norm distribution $\NNN_y$ 
(related to $F$) is called the \emph{universal Euler system} of level $(F,y)$.
This is the concept introduced by Rubin in his book Euler 
systems~\cite{Rubin2}.

\subsection{Function field case} Let $K=F_q(T)$ and $A=F_q[T]$. For any
$f(T)\in A$, let $K(f)=K(\lambda_f)$ be the cyclotomic function field 
of $K$ related to $f$ where $\lambda_f$ is a division point of $f$
with respect to the Carlitz module. The Galois group $G_f$ of
$K(f)/K$ is known to be isomorphic to $(A/f)^{\times}$. Thus we can 
identify every $\sigma=\sigma_x\in G_f$ for some(a unique) 
$x\in  (A/f)^{\times}$.
The \emph{ordinary distribution of level $f$} on the function 
field $K$ is defined to be a map
\[ \phi: \frac{1}{f}A/A\longrightarrow Ab=\text{abelian group} \]
satisfying
\[ \phi(x)=\sum_{py =x} \phi(y), \forall  p\mid f, x\in \frac{p}{f}A/A. \]
One can then talk about the \emph{universal ordinary distribution} as
the universal object to the category of ordinary distributions. As in the
number theory counterpart, by abusing notation, we say the group
\[ U_f= \frac{\langle [a]: a\in \frac{1}{f}A/A \rangle}
{\langle [a]-\sum_{pb=a}[b]: p\mid f, a\in \frac{p}{f}A/A\rangle} \]
the universal ordinary distribution. $U_f$ is naturally equipped with
a $G_f$-action by sending $\sigma_x[a]=[xa]$.
This distribution plays a similar
role to the universal ordinary distribution in the study of cyclotomic
function field. One can easily show that $U_f$ is a free abelian group
of order $|G_f|$.
We'll see that it actually is also a special case 
of the universal norm distribution. 

We work on more generality. Let $K$ be a fixed function field.
Pick a place $\infty$ in $K$. Let $A$ be the integer ring
corresponding to the place $\infty$. Choose a sign function $\sgn$ 
on $K^{\ast}$. Let $\phi$ be a sign-normalized Drinfeld module of rank $1$.
The field $H^{+}$ is defined to be the extension of $K$ by adding all
the coefficients of $\phi_a$ for $a\in A$.

For any ideal $I$ of $A$, 
let $K(I)$ be the cyclotomic function field extension of $K$ related 
to $I$(and related to the sign-normalized Drinfeld module $\phi$). 
Let now $X$ be the set of all prime ideals of $A$, then $Z$ can be
considered as the set of all integral ideals of $A$. Let $G_I$ be the
Galois group of $K(I)/H^{+}$. We know that $G_I=(A/I)^{\times}$ and thus 
$G_I$ satisfies the condition of the universal norm distribution. For any
$\wp\in X$, we can define the Frobenius element $\Fr_{\wp}$ correspondingly.
Let $\co=\co_{H^{+}}$ and $\Phi=H^{+}$. We can now define the universal
norm distribution by choosing a set of polynomials $\{p(\wp,t)\}$. 

In particular, if let $K=F_q(T)$ and let $A=F_q[T]$. Let the sign normalized
Drinfeld module be the usual Carlitz module. In this case $H^{+}$ is 
actually just $K$. let  $p(\wp,t)=1-t$ for every $\wp\in X$. 
We can identify the 
corresponding universal norm distribution $\NNN_f$ related to $f\in A$
as the $G_f$-module
\[ \NNN_f= \frac{\langle [\sigma f']: f'\mid_s f, \sigma\in G_{f'}\rangle}
{\langle (1-\Fr^{-1}_p)[\sigma f']-N_{f(p)}[\sigma f(p)f']: f(p)f'\mid_s f, 
\sigma\in G_{f'}\rangle}. \]
Now we can define a homomorphism from $\NNN_f$ to $U_f$ by sending
$[f']\longmapsto [\frac{1}{f}]$. This homomorphism is shown to be 
an isomorphism.

\section{Basic properties of the universal norm distribution $\NNN_Z$}
\label{section:basic}

Recall by our definition, for every $z\in Z$,
$\AAA_z$ is a free $\T$-module generated by the set
\[ A_z=\bigcup_{\substack{z'\ \text{finite}\\ z'\mid_s z}} B_{z'}=
\bigcup_{\substack{z'\ \text{finite}\\ z'\mid_s z}} \{[gz']:g\in G_{z'}\}. \]
If let $B_n$ be the set of all elements  
\[ \{[g z]\in A: \text{the restriction $g_{z(x)}=1$ 
for exactly $n$ primes $x\mid z'$}\} \] 
Then $A_z$ is the disjoint union   
\[ A_z=\bigcup_{n\geq 0}\bigcup_{\substack{z'\ \text{finite}\\z'\mid_s z} } 
\left ( B_n\cap B_{z'} \right ). \]
We have the following key proposition:

\begin{Proposition} \label{Proposition:basis}
The free $\T$-module $\AAA_z$, for every $z\in Z$,
possesses a $\T$-basis 
\[ \{ \lambda_{z''} [g z']: z', z' z''\mid_s z,\ z'z''\in \zf,\
[g z']\in B_0 \} \]
where $\lambda_{z''}$ is defined to be the product of $\lambda_{z(x)}$ for
all $x\mid z''$.
\end{Proposition}
\begin{proof} Suppose that $[g z']\in B_n\cap A_z$ for
$n\geq 1$, then there exists a prime $x\mid z'$ such that 
$g_{z(x)}=1$. One has
\[ [g z']=-\sum_{1\neq g'\in G_{z(x)}} [g g' z']
-\lambda_{z(x)} [g z'/z(x)]+p(x; \Fr^{-1}_x) [g z'/z(x)]. \]
Thus
\[ \langle B_n\rangle_{\T} \cap \AAA_z\subseteq
\langle B_{n-1} \rangle_T \cap \AAA_z+\sum_{x\mid z} 
\lambda_{z(x)} A_{z/z(x)}+\sum_{x\mid z} \AAA_{z/z(x)} \]
where $\langle B_n \rangle_{\T}$ denotes the free $\T$-module generated 
by $B_n$. Thus by induction, the set given in the proposition 
generates $\AAA_z$. 
We just need to show the cardinality of this set agrees with the $\T$-rank
of $\AAA_z$. For finite $z\in Z$, the $\T$-rank of $\AAA_z$ is 
\[ \sum_{z'\mid_s z} |G_{z'}|=\prod_{x\mid z}(|G_{z(x)}|+1). \]
On the other hand, the cardinality of the set in the proposition is
\[ \begin{split}
 \sum_{z''\mid_s z} \sum_{z'\mid_s \frac{z}{z''}} |B_0\cap B_{z'}|
=& \sum_{z''\mid_s z }
\sum_{z'\mid_s \frac{z}{z''}} \prod_{x\mid z'} (|G_{z(x)}|-1)\\
=& \sum_{z''\mid_s z} 
\prod_{x\mid \frac{z}{z''}} |G_{z(x)}| \\
=& \prod_{x\mid z}(|G_{z(x)}|+1).
\end{split} \]
This proved the case when $z$ is finite. Taking the limit, then we have
the proof for infinite $z\in Z$.
\end{proof}
\begin{Proposition} \label{Proposition:basis2}
(1). The module $\NNN_z$ is a free $\T$-module with basis
$B_0\cap A_z$.

(2). For every $z'\mid_s z$, the natural injection of $\NNN_{z'}$ to
$\NNN_z$ is a splitting $G_z$-monomorphism.

\end{Proposition}
\begin{proof} Immediately from  Proposition~\ref{Proposition:basis}.
\end{proof}
\begin{rem} From the above Proposition~\ref{Proposition:basis2}(1), one see 
that $\NNN_z$ is free
$\T$-module of rank $|G_{z}|$. In particular, in the trivial universal norm 
distribution case, one
see that the image of $B_z$ in $\NNN_z$ actually is a basis of $\NNN_z$,
thus $\NNN_z$ is isomorphic to $\T[G_z]$, which justifies the meaning of
\emph{trivial}.
\end{rem}

\begin{rem} \label{rem:basis}
From the above Proposition~\ref{Proposition:basis2}(2), we'll
henceforth identify $\NNN_{z'}$ as a submodule of $\NNN_z$. In particular, for
every $z\in Z$, we have
\[ \NNN_{z}=\bigcup_{\substack{z'\ \text{finite}\\ z'\mid_s z}} \NNN_{z'}. \] 
This observation will be used to the study of the universal Kolyvagin 
recursion in Ouyang~\cite{Ouyang4}.
\end{rem}

\begin{Proposition} \label{Proposition:embed}
Let $w\mid z$ be a pair of elements in $Z$. Then the
corestriction homomorphism $\cor_{w,z}$ from $\AAA_w$ to $\AAA_z$ by
\[ [w']\longmapsto N^{z'}_{w'}[z'](w'\mid_s w, z'\mid_s z, \bar w'=\bar z') \]
induces an embedding from $\NNN_w$ to $\NNN_z$. 
In particular, when
$w\mid_s z$, this embedding is the natural injection as given in 
Proposition~\ref{Proposition:basis2}.
\end{Proposition}
\begin{proof}
Write $V_1$(resp. $W_1$) the free $\T$-submodule of $\AAA_w$(resp. $\AAA_z$) 
generated by $B_0\cap A_w$(resp. $B_0\cap A_z$). Write $V_2$(resp. $W_2$)
the free  $\T$-submodule of $\AAA_w$(resp. $\AAA_z$) generated by other 
elements in the basis of $\AAA_w$(resp. $\AAA_z$) given 
by Proposition~\ref{Proposition:basis}. Then it is easy to check that
$\cor_{w,z}$ maps $V_i$ to $W_i$ injectively. Hence it induces a well 
defined embedding from $\NNN_w$ to $\NNN_z$.
\end{proof}

\section{Anderson's resolution} \label{section:and}

\subsection{Set up} 
Let $z\in Z$ be given. Let
\[ \LLL_z=\bigoplus_{y\mid \bar z}
\AAA_{z/z(y)}[y] \]
where $y$ is finite and $[y]$ is a symbol depending only on $y$.  
If we write
\[ [g'z'][y]=[g'z', y]  \]
for elements in $\AAA_{z/z(y)}[y]$, 
then $\LLL_z$ is the  free $\T$-module generated by the set
\[ \{ [a,y]:  [a]\in A_{z/z(y)}, y\mid \bar z \} \]
We assign a grade in $\LLL_z$ by declaring 
\[ \deg [a,y]= -\deg  y. \]
For any
$g\in G_z$ and $[g'z']\in A_{z/z(y)}$, declare the $G_z$-action as
\[ g[g'z',y]:= [g_{z'} g'z', y], \]
then $\LLL_z$ becomes a graded $\T[G_z]$-module. 
$\LLL_z$ is bounded above since all its non-negative components are $0$.
Moreover,  $\LLL_z$ is bounded if and only if $z$ is finite.

With abuse of notation,
denote by  $\lambda_{z(x)}$, $\lambda_{z'}$ the homomorphisms of $\LLL_z$ 
inheriting from the homomorphisms in $\AAA_z$ bearing the same names. 
Now let
\[ d: \LLL_z \longrightarrow \LLL_z,\ [a,y]\longmapsto
\sum_{x\mid y} \omega(x,y)\lambda_{z(x)}[a, y/x] \]
where $\omega$ is as defined in \S~\ref{subsection:Notations}.
Clearly $d$ commutes with $G_z$-actions.
A straightforward calculation shows that $d^2=0$ and therefore $d$ is a 
differential of degree $1$.  Define an $\T[G_z]$-homomorphism
$\mathbf u: \LLL_z\rightarrow \NNN_z$ by
\[ [a,y]\longmapsto \begin{cases} [a],\ &\text{if}\ y=\mathbf 1;\\
0,\ &\text{if}\ y\neq \mathbf 1. \end{cases} \]
Regard $\LLL_z$ as a complex $\LLL^{\bullet}_z$ by the differential $d$, 
and regard $\NNN_z$ as a complex concentrated on $0$-component.
Then one can easily check that $\mathbf u$ is a homomorphism of complexes.
Because of the following Theorem,
we call the complex $(\LLL^{\bullet}_z, d)$
(or simply $\LLL^{\bullet}_z$)
\emph{Anderson's resolution} of the universal norm system $\NNN_z$. 

\begin{theo} \label{Proposition:Anderson}
The homomorphism $\mathbf u$ is a quasi-isomorphism, i.e., 
the complex $(\LLL^{\bullet}_z,d)$ is acyclic for degree
$n\neq 0$ and   $H^0(\LLL^{\bullet}_z, d)\cong\NNN_z$ induced by $\mathbf u$.
\end{theo} 
\begin{proof} For any $a\in B_0\cap B_{z/z(y)}$, 
consider the graded $\T$-submodule $C^{\bullet}_a$ of $\LLL^{\bullet}_z$ 
generated by
\[ \{ \lambda_{w}[a,y'],\ w\mid_s z,  \bar w y' \mid y \}. \]
One can see that $C^{\bullet}_a$ is $d$-stable. Thus $C^{\bullet}_a$ is 
actually a subcomplex of $\LLL^{\bullet}_z$. By 
Proposition~\ref{Proposition:basis}, 
$\LLL^{\bullet}_z$ is the direct sum of
$C^{\bullet}_a$ for $a$ over $B_0\cap A_z$. We hence only have to study the 
complex  $C^{\bullet}_a$.
Now the theorem follows from Lemma~\ref{Lemma:iso}.
\end{proof}

\subsection{The Koszul complex $\tilde C^{\bullet}_{y}$} 
Let $\Lambda$ be the polynomial ring
\[ \Lambda=\T[Z]=\{\sum t_z z: t_z\in T, z\in Z\}.\]
Let $\tilde C^{\bullet}_{y}$ be
the Koszul complex of $\Lambda$ with the regular sequence
$x_1<\cdots<x_m$ where $y=x_1\cdots x_m$. 
Thus $\tilde C^{\bullet}_{y}$ is the graded exterior algebra
\[ \bigoplus_{y'\mid y} 
\Lambda e_{y'} \]
with 
\[ e_{y'}=e_{x_{i_1}}\wedge \cdots \wedge e_{x_{i_k}},\ \text{and}\ 
\deg e_{y'}=-\deg y'=-k \]
where 
\[ y'= x_{i_1}\cdots x_{i_k}, x_{i_1}<\cdots< x_{i_k}.\]
The corresponding  differential is given by
\[ d\, e_{x}=x. \]

\subsection{Truncated Koszul subcomplex $C^{\bullet}_{y}$}
Let $C^{\bullet}_{y}$ be the graded $\T$-submodule of
$\tilde C^{\bullet}_{y}$ generated by all elements of the
form $y'' e_{y'}$ for all $y'y''\mid y$. This submodule
is stable under the differential, thus is a subcomplex of
$\tilde C^{\bullet}_{y}$. Moreover, it is a direct summand of 
$\tilde C^{\bullet}_{y}$. By the general theory of Koszul complex,
$C^{\bullet}_{y}$ is acyclic in nonzero degree and 
$H^0(C^{\bullet}_{y})$ is a free $\T$-module generated by 
$e_{\mathbf 1}$.

\begin{Lemma} \label{Lemma:iso}
For any $a\in B_0\cap B_{z/z(y)}$, the complex
$C^{\bullet}_a$ is isomorphic to $C^{\bullet}_{y}$. Thus
$C^{\bullet}_{a}$ is acyclic in nonzero degree and 
$H^0(C^{\bullet}_{a})$ is a free $\T$-module generated by 
$[a,1]$.
\end{Lemma}
\begin{proof} Let $C^{\bullet}_{y}$ act on $C^{\bullet}_{a}$ by
\[ x [a,y']=\lambda_{z(x)}[a,y'] \]
and
\[ e_{x}[a,y']=\begin{cases} (-1)^{|\{x'<x: x'\mid y'\}|} 
[a, x y']\ &\text{if}\ x\nmid y';\\
0\ &\text{if}\ x\mid w. \end{cases}  \]
By straightforward calculation  
\[ d(\xi\eta)=(d\xi)\eta+(-1)^{\deg \xi}\xi(d\eta), 
\xi\in C^{\bullet}_{y}, \eta\in C_x^{a}. \]
Thus $C^{\bullet}_a=C^{\bullet}_{y} [a,1]$.
\end{proof}

\subsection{Compatibility } \label{subsection:embed}
From Proposition~\ref{Proposition:embed}, the injective
corestriction homomorphism $\cor_{w,z}$ from $\AAA_w$ to $\AAA_z$
induces a corestriction homomorphism $\cor$ 
from $\LLL_w$ to $\LLL_z$ by
\[ [a,y]\longmapsto \cor_{w/w(y),z/z(y)}[a, y]. \]
A straightforward calculation shows that $\cor$ is 
compatible with the differential $d$.
Now if let $\tilde \LLL_z$ be the extended exact sequence of $\NNN_z$ to 
$\LLL_z$, i.e., $\tilde \LLL_z$ is the sequence
\[ \cdots \LLL^{-n}_z\rightarrow \cdots \rightarrow \LLL^0_z
\xrightarrow{\mathbf u} \NNN_z \rightarrow 0 \]
then the corestriction map $\cor$ is actually an injective chain homomorphism
from $\tilde \LLL_w$ to $\tilde \LLL_z$ and is thus an embedding.
When $w\mid_s z$, this embedding
$\cor$ is again a natural injection.

\subsection{Connecting map for different norm distributions} Now fix $X$ 
and $\T$, suppose that we have two sets 
of polynomials $\{p_1(x;t)\}$ and $\{p_2(x;t)\}$ in $\co[t]$, 
then we have two norm distributions $\NNN_{1,z}$ and  $\NNN_{2,z}$,
and two corresponding Anderson's
resolutions $\LLL_{1,z}$ and  $\LLL_{2,z}$.
Then there exists a connecting homomorphism
\[ \phi_{1,2}: \LLL_{1,z}\otimes_{\co} \Phi \longrightarrow 
\LLL_{2,z}\otimes_{\co} \Phi \]
by
\[ [z',y]\longmapsto \sum_{w\mid_s z'}(-1)^{\deg \bar w} \prod_{x\mid w}
\frac{p_2(x,\Fr^{-1}_x)-p_1(x,\Fr^{-1}_x)}{|G_{z(x)}|}[z'/w,y]. \]
By straightforward calculation, one can
check that $\phi_{2,1}$ is the inverse of $\phi_{1,2}$, thus
$\phi_{1,2}$ is actually an isomorphism, which induces isomorphisms between
$\NNN_{1,z}\otimes_{\co} \Phi$ and $\NNN_{2,z}\otimes_{\co} \Phi$.
In particular, if we let 
$p_1(x;t)\equiv 1$ for every $x\in X$, then $\NNN_{1,z}\otimes_{\co} \Phi$
is nothing but the module $\T[G_z]$, thus we have
\begin{Proposition}
The $\T\otimes_{\co}\Phi[G_z]$ module $\NNN_{z}\otimes_{\co} \Phi$ 
is free of rank $1$ for every universal norm distribution.
\end{Proposition}

\subsection{Double complex structure of $\LLL_z$} \label{subsection:doubleL}
Set a bidegree in $\LLL_z$ by
\[ \deg^{(2)} [z',y]=(\deg \bar z', -\deg \bar z'-\deg y). \]
We set 
\[ \begin{split}
& d_{1,x}[ z',y]= -\omega(x,y)
N_{z(x)}[z'z(x),y/x],\ \\ 
& d_{2,x}[z',y]= \omega(x,y)
p(x;\Fr_x)[z', y/x]. \end{split} \]
and let 
\[ d_{x}=d_{1,x}+d_{2,x},\ \ d_{1}=\sum_x d_{1,x}, \ 
d_2= \sum_x d_{2,x}. \]

\begin{Lemma} (1). For any $x, x'\mid z$, $i=1, 2$, 
\[ d^2_{i,x}=d_{1,x} d_{2,x'}+ d_{2,x'} d_{1,x}=0. \]

(2). $d_1^2=d^2_2=d_1d_2+d_1 d_2=0$.

(3). $d_{i,x}$ is $G_z$-stable.
\end{Lemma}
\begin{proof} Straightforward. \end{proof}
From the above lemma, we see that $\LLL_z$ is equipped with 
a multiple complex structure. 
In particular, $(\LLL^{\bullet,\bullet}_z;d_1,d_2)$ is a double complex 
corresponding to the above bigrading. 
We'll use this complex to study the group cohomology of $\NNN_z$ 
in \S~\ref{sec:ggt}.

\section{Preparation from homological algebra}
\label{sec:prep}

\subsection{Complex of type $E$} 
\label{subsection:exterior}
Let $A$ be a free $\co$-module of finite rank.
Let $\Lambda_A=\Lambda_A(x_1,\cdots, x_t)$ be the exterior algebra over $A$, 
with the differential $d$ given by $d(x)= \sum_i m_i x\wedge x_i$ where 
$m_i\in \co$.  For each $S\subseteq\{1,\cdots,t\}$, let $m_S$ be the 
greatest common divisor of $m_i$ for all $i\in S$. In particular, let
$m$ be the greatest common divisor of $m_i$ for all $1\leq i\leq t$.

Let $S=\{i_1,\cdots, i_s\}$ such that $i_1\leq \cdots \leq i_s$.
Let $\{e_S=x_{i_1}\wedge\cdots \wedge x_{i_s}\}$ be the standard basis 
of $\Lambda_A$. By linear algebra, in the $\Phi$-vector space generated by 
$\{x_1,\cdots, x_t\}$, there exists another basis $\{y_1,\cdots,y_t\}$
such that $y_1=\frac{1}{m}\sum_i m_i x_i$ and the transformation matrix is
inside $SL(t,\ZZ)$, thus $\{ e'_S=y_{i_1}\wedge \cdots\wedge y_{i_s}\}$ is 
another basis for $\Lambda_A$. Hence one can easily show that 
$H^{\ast}(\Lambda_A)$ is a free graded $A/mA$-module generated by cocycles
represented by $e'_S$ for all $S$ which contains $1$, thus is 
a free $A/mA$-module of rank $2^{t-1}$, with the $i$-th component
a free $A/mA$-module of rank $\binom{t-1}{i-1}$(or $0$ if $i=0$).

\subsection{The tensor projective resolution $P_{z\bullet}$} This setup
is from Ouyang~\cite{Ouyang2}. Fix an element $z\in Z$.
Assuming that $G_{z(x)}$ is a cyclic group for every $x\mid z$. Let 
$\sigma_{z(x)}$ be a generator of $G_{z(x)}$. It is well known that
the sequence
\[ \cdots \ZZ[G_{z(x)}]\xrightarrow{N_{z(x)}}
\ZZ[G_{z(x)}]\xrightarrow{1-\sigma_{z(x)}}\ZZ[G_{z(x)}]
\xrightarrow{\epsilon}\ZZ\rightarrow 0\]
is exact, where $\epsilon$ is the augmentation map. Let $P_{z(x)\bullet}$
be the resulting resolution for the trivial $\ZZ[G_{z(x)}]$-module $\ZZ$, 
we can thus write $P_{z(x)\bullet}$ as the graded module
\[ \bigoplus_{n\geq 0}  \ZZ[G_{z(x)}][x^n] \]
with the symbol $[x^n]$ is of degree $n$ and the differential given by
\[ \partial_{z(x)}[x^n]=\begin{cases} 
(1-\sigma_{z(x)})[x^{n-1}],&\ \text{if}\ n>0\ \text{odd};\\
N_{z(x)}[x^{n-1}],&\ \text{if}\ n>0\ \text{even}. \end{cases} \]
Now let $P_{z\bullet}$ as the tensor
product of $P_{z(x)\bullet}$ over all $x\mid z$. $P_{z\bullet}$ is 
the so called \emph{tensor projective resolution} of the trivial 
$\ZZ[G_z]$-module $\ZZ$ with respect to the cyclic decomposition
\[ G_z=\prod_{x\mid z} G_{z(x)}=\prod_{x\mid z} \langle \sigma_{z(x)} 
\rangle. \]
Let $[w]$ be an indeterminate for every $w\in Z$. 
Then the tensor resolution $P_{z\bullet}$
is the projective 
$\ZZ[G_z]$-resolution of the trivial module $\ZZ$ by
\[ P_{z, n}=\bigoplus_{\substack{\bar w\mid z\\ \deg w=n}}
\ZZ[G_z] [w] \]
and the differential $\partial_z$
is given by
\[ \partial_z [w]=\sum_{x\mid w} 
(-1)^{\sum_{x'<x}v_{x'}w} \alpha_{z(x)}[w/x] \]
where $\alpha_{z(x)}$ is equal to $\sigma_{z(x)}-1$ if $v_{x} w$ odd
and $N_{z(x)}$ if $v_{x} w$ even. For any $z'\mid_s z$,
one has a natural inclusion of $P_{z'\bullet}$ to $P_{z\bullet}$
by sending $[w]$ to $[w]$.

\subsection{$G_z$-cohomology of trivial module $A$} Let $A$ be 
a free $\co$-module with trivial $G_z$-structure.
To compute its $G_z$-cohomology, it suffices to
compute the cohomology
\[ I^{\bullet}_{A, z}=\Hom_{\ZZ[G_z]}(P_{z\bullet},A)=
\bigoplus_{\substack{w\ \text{finite}\\ \bar w\mid z}} A[w] \]
with the differential
\[ \delta_z[w]=\sum_{x\mid z} (-1)^{\sum_{x'<x}v_{x'}w}
a_{z(x)} [wx] \]
where $a_{z(x)}$ is equal to $0$ if $v_{x}w$ even and to
$|G_{z(x)}|$ if $v_{x}w$ odd. The inclusion of $P_{z'\bullet}$ 
to $P_{z\bullet}$ for $z'\mid_s z$ thus induces a projection from 
$I^{\bullet}_{A, z}$ to $I^{\bullet}_{A, z'}$. One see that 
$I^{\bullet}_{A, z'}$ is a direct 
summand of $I^{\bullet}_{A, z}$. 

For any finite $w$ with $\bar w\mid z$, let
\[ I^{\bullet}_A [w^2]=\bigoplus_{w'\mid \bar w} A[w^2/w'], \]
then $I^{\bullet}_A[w^2]$ is a direct summand of $I^{\bullet}_{A, z}$
and 
\[ I^{\bullet}_{A, z}=\bigoplus_{\bar w\mid \bar z} I^{\bullet}_A[w^2]. \]
If $w=1$, the subcomplex $ I^{\bullet}_A[w^2]$ is just a copy of
$A$ with the differential $0$, thus the cohomology of it is $A$ too.
If $w\neq 1$, the subcomplex $ I^{\bullet}_A[w^2]$ is of type $E$. 
Let $m_w$ be the
greatest common divisor of $|G_{z(x)}|$ for $x\mid w$, then 
$H^{\ast}( I^{\bullet}_A[w^2])$ is then a free graded $A/m_w A$-module
of rank $2^{\deg \bar w-1}$. One see the $(2\deg w-\deg \bar w+i)$-th
cohomology is just a free $A/m_w A$-module of rank
$\binom{\deg \bar w-1}{i-1}$ for $1\leq i\leq \deg w$ and $0$ 
otherwise.

Denote $H^{\ast}( I^{\bullet}_A[w^2])$ by $H_{A,w}$. Then
with the above analysis, one  has

\begin{Proposition}  \label{Proposition:trivial1} Fix a finite $z\in Z$
such that every $G_{z(x)}$ is cyclic for $x\mid z$. For a free
$\co$-module $A$ with trivial $G_z$-action, then we have

(1). For any $z'\mid_s z$, the cohomology group
$H^{\ast}(G_{z'}, A)$ is a direct summand of $H^{\ast}(G_z,A)$.

(2). The cohomology group $H^{\ast}(G_z,A)$ is the direct sum of 
$H_{A,w}$ for every $\bar w\mid \bar z$ where: (a). For $w=1$, $H_{A,w}=A$
is with grade $0$; (b). For $w\neq 1$, $H_{A,w}$ is a free graded 
$A/m_w A$-module with the $(2\deg w-\deg \bar w+i)$-th component
of rank $\binom{\deg \bar w-1}{i-1}$ for $1\leq i\leq \deg \bar w$
and $0$ for otherwise.

\end{Proposition}
\begin{rem}

\end{rem}
Now for a finite fixed $z\in Z$, suppose $M\in \co$ a common divisor 
of $|G_{z(x)}|$ for every $x\mid z$. Then
the case for $G_z$-cohomology of $A/MA$ is much simpler. In this case,
\[ H^{\ast}(G_z, A/MA)=H^{\ast}(I^{\bullet}_{A,z}/M I^{\bullet}_{A,z}), \]
and the differential in $I^{\bullet}_{A,z}/M I^{\bullet}_{A,z}$ is nothing
but $0$, thus $H^{\ast}(G_z, A/MA)$ as a graded module is isomorphic
to $I^{\bullet}_{A,z}/M I^{\bullet}_{A,z}$.
One  has

\begin{Proposition} \label{Proposition:trivial}
There exists a family
\[ \{ [w]\in H^{\ast}(G_z, A/MA): w\ \text{finite},\bar w\mid z \} \]
with the following properties:

(1). For any $z'\mid_s z$, the restriction of the family
\[ \{ [w]: \bar w\mid \bar z',\deg w=n\} \]
to $H^n(G_{z'}, A/MA)$ is an $A/MA$-basis of the latter one.

(2). The restriction of $[w]$ for $\bar w\nmid \bar z'$ to 
$H^{\ast}(G_{z'}, A/MA)$ is $0$.
\end{Proposition}

\section{$G_z$-cohomology of the universal norm distribution $\NNN_z$}
\label{sec:ggt}
In this section, we use tools developed in the previous sections
to study the $G_z$-cohomology of the universal norm distribution $\NNN_z$
and of $\NNN_z/ M\NNN_z$. We assume that $G_{z(x)}$ cyclic for every 
$z\in Z$ and $M$ a common divisor
of $|G_{z(x)}|$ for every $x\mid z$.

\subsection{Setup of double complex $\K^{\bullet,\bullet}_{z}$} 
With preparations from the above two sections, we let
\[ \K^{\bullet,\bullet}_{z}=\Hom_{G_z}(P_{z\bullet},\LLL^{\bullet}_z)
\]
If we write $[a,y,w]=([w]\mapsto [a,y])$, then 
$\K^{\bullet,\bullet}_{z}$ 
is the free graded $\T$-module with basis
\[ \{ [a,y,w]: y\mid \bar z, a\in A_{z/z(y)}, \bar w\mid \bar z \} \]
and with the double grading given by 
\[ \deg[a,y,w]=(-\deg y, \deg w). \]
The induced $\T[G_z]$-module structure is given by 
\[ g[a, y, w]=[gx, y, w] \]
for any $g\in G_z$. Use the same notations for the
operators in $\K^{\bullet,\bullet}_{z}$ induced from 
$\LLL^{\bullet}_z$, i.e., $\lambda_{z(x)}$,
$\lambda_{z}$ and so on.
Now the two differentials of $K^{\bullet,bullet}_z$ are given by
\[ d[a,y,w]=\sum_{x\mid y}\omega(x,y) \left (
p(x;\Fr^{-1}_x) [a,y/x,w]-N_{z(x)}[z(x)a,y/x,w] \right ), \]
\[ \delta[a, y, w]=(-1)^{\deg y}\sum_{x\mid z}
(-1)^{\sum_{x'<x}v_{x'}(w)} a_{z(x)}[a, y, wx] \]
where $a_{z(x)}$ is equal to $1-\sigma_{z(x)}$ if $v_{x}(w)$ even
and $N_{z(x)}$ if $v_{x}(w)$ odd.
Let $\K^{\bullet}_{z}$ be 
the single total complex of $\K^{\bullet,\bullet}_{z}$.
and let $\K_{z}$ be the underlying module.

Let $\bar \K^{\bullet}_{z}=
\Hom_{G_z}(P_{z\bullet},\NNN_z)$. Then
it is the quotient of free $\T$-module generated by
\[ \{[a,w], a\in A_z,\bar w\mid \bar z \} \]
modulo relations generated by
\[ \lambda_{z(x)}[a,w],\ a\in A_{z/z(x)},\ \bar w\mid \bar z, 
\forall\ x\mid z, \]
with the differential $\delta$ given by
\[ \delta [a,w]=\sum_{x\in z} (-1)^{\sum_{x'<x}v_{x'}w}
a_{z(x)} [wx]. \]
We have the induced map
\[ \mathbf u: \K^{\bullet}_z\longrightarrow
\bar \K^{\bullet}_z,\  [a,y,w]\longmapsto
\begin{cases} [a,w],\ &\text{if}\ y=\mathbf 1; \\
0,\ &\text{if}\ y\neq \mathbf 1. \end{cases} \]

\begin{Proposition} \label{Proposition:cohomology}
The homomorphism $\mathbf u$ is a quasi-isomorphism. Thus

(1). $H^{\ast}(\K^{\bullet}_{z}, d+\delta)\cong 
H^{\ast}(G_z, \NNN_z)$.

(2). $H^{\ast}(\K^{\bullet}_{z}/M\K^{\bullet}_{z}, 
d+\delta)\cong H^{\ast}(G_z, \NNN_z/M\NNN_z)$.
\end{Proposition}
\begin{proof}
By Theorem~\ref{Proposition:Anderson}, $\ker \mathbf u$ is $d$-acyclic,
by spectral sequence argument, it is hence $(d+\delta)$-acyclic. Thus 
$\mathbf u$ is a quasi-isomorphism. (1) follows immediately from
the quasi-isomorphism. Since both $\K^{\bullet}_{z}$ and
$\NNN_z$ are free $\T$-modules, the induced homomorphism
$\bar {\mathbf u}$ from $\K^{\bullet}_{z}/M\K^{\bullet}_{z}$
to $\bar \K^{\bullet}_{z}/M\bar \K^{\bullet}_{z}$
is also a quasi-isomorphism and (2) follows immediately.
\end{proof}

\subsection{Another double complex structure of $\K_z$}
\label{subsection:another}
Keep $\K_z$ as the same bigraded module as in the previous section. 
Let's equip it with different differentials 
$(\tilde d,\tilde \delta)$ as the following:
\[ \tilde d[a,y,w]=\sum_{x\mid y}\omega(x,y) (-1)^{\sum_{x'<x}v_{x'}(w)}
\left (p(x;\Fr^{-}_x) [a,y/x,w]-N_{z(x)}[az(x),y/x,w] \right ), \]
\[ \tilde \delta[a, y, w]=\sum_{x\mid z} (-1)^{\sum_{x'\leq x}v_{x'}(y)}
(-1)^{\sum_{x'<x}v_{x'}(w)} a_{z(x)}[a, y, wx].\]
One can easily check that 
\[ \tilde d^2=\tilde \delta^2=\tilde d\tilde \delta+\tilde \delta\tilde d=0. \]
We define an involutive $G_z$-equivariant bigraded automorphism $\epsilon$
of $\K_z$ by the rule
\[ \epsilon[a,y,w]=(-1)^{\sum_{x,x': x'<x} v_x(y)v_{x'}}(w). \]
by a straightforward calculation, one finds that
\[ \epsilon\ \tilde d\ \epsilon=d, \qquad 
\epsilon\ \tilde \delta\ \epsilon=\delta. \]
Thus $\epsilon$ induces an isomorphism between 
the cohomology of $(\K;\tilde d,\tilde \delta)$ and 
the cohomology of $(\K;d,\delta)$, which is then isomorphic to the 
$G_z$-cohomology of $\NNN_z$. 

In the sequel, we'll use the double complex $(\K;d,\delta)$ to study the
cohomology of $\NNN_z$. However, the results obtained here is easy to adapt
to the double complex $(\K;\tilde d,\tilde \delta)$.
The double complex $(\K;\tilde d,\tilde \delta)$
will be used to the study of the universal 
Kolyvagin recursion in Ouyang~\cite{Ouyang4}.

\subsection{Multiple complex structure of $\K_z$} The underlying
module $\K_z$ has abundant complex structures. For $x\mid z$, set
\[ \deg_{1,x}([z', y, w]):=v_{x} (z'), \]
\[ \deg_{2,x}([z', y, w]):=v_{x}(yz'), \] 
\[ \deg_{3,x}([z', y, w]):=v_{x}(w). \]
We call $\deg_{i,x}([z', y, w])$ for $i=1,2,3$ the
$(i,x)$-degree of $[z', y, w]$. Make the degrees invariable with
$G_z$ action, then
$\K_{z}$ is equipped with a multi-graded module structure.  Let
\[ d_{1,x}[a, y, w]:= -\omega(x,y) 
N_{z(x)}[az(x), y/x, w] \]
\[ d_{2,x}[a,y,w]:= \omega(x,y) 
p(x;\Fr^{-1}_x) [a, y/x,w] \]
\[ d_{3,x}[a,y,w]:=(-1)^{\deg y}
(-1)^{\sum_{x'<x}v_{x'}w} a_{z(x)}[a, y, wx]. \]
The map $d_{i,x}$ is of $(i,x)$-degree $+1$. 
It is easy to check that for every $i$, $j=1,2,3$ and $(i,x)\neq (j,x')$, 
one has
\[ d^2{i,x}=d_{i,x}d_{j,x'}+d_{j,x'}d_{i,x}=0. \]
Thus $d_{i,x}$ are differentials 
of $\K_z$ observing the above multi-grading structure.
One see that $d$ is the sum of all $d_{i,x}$ for $i=1,2$ and
$x\mid z$ and $\delta$ is the sum of $d_{3,x}$. The total degree of
$\K_z$ is just the sum of all $(i,x)$-degrees. Thus we can use
this multi-complex structure to study the total cohomology
of $\K_z$ and hence the $G_z$-cohomology of $\NNN_z$.

Furthermore, note that any combination of $d_{i,x}$ is still 
a differential in $\K_z$. In particular, 
$d_{i}=\sum_{x\mid z} d_{i,x}$ for $i=1, 2$
is the differential induced by the differential $d_i$ in
$\LLL_z$ when viewing $\LLL_z$ as a double complex. We have
$d=d_1+d_2$ and $\delta=\sum_{x\mid z} d_{3,x}$. Correspondingly,
we can make $\K_z$ as a triple complex $\Kc_z$ with
differentials $d_1$, $d_2$ and $\delta$. As a convention, we use
$m, n, p=m+n$ and $q$ to denote the corresponding degrees for
the differentials $d_1, d_2, d$ and $\delta$.
We shall use this triple complex structure
of $\K_z$ to study the total
cohomology of $\K^{\bullet}_z$.

\subsection{Compatibility} For every $z'\mid_s z$, let $\K_{z'}$ be the 
submodule of $\K_{z}$ generated by 
\[ \{[a, y, w]: y\mid z', a\in B_{z'/z(y)}, 
\bar w\mid z'\} \]
and let $\K_{z}(z')$ be the submodule generated by
\[ \{[a,y,w]:y\mid z', a\in B_{z'/z(y)}, 
\bar w\mid z\} \]
One can check that $\K_{z'}$ and $\K_{z}(z')$ are compatible
with differentials. The $(d+\delta)$-cohomology of $\K_{z'}$ is 
just $H^{\ast}(G_{z'}, \NNN_{z'})$ and the $(d+\delta)$-cohomology
of $\K_{z}(z')$ is $H^{\ast}(G_z, \NNN_{z'})$. Moreover, if using the 
embedding $\cor$ defined in \S~\ref{subsection:embed} for Anderson's 
resolution, then for every 
$w\mid z$, one has a well defined embedding from $\K_{w}$ to $\K_{z}$.

\subsection{The study of spectral sequences}
We now discuss the $G_z$-cohomology of $\NNN_z$ and 
$\NNN_z/M\NNN_z$. 
We study the triple complex 
$(\K^{\bullet,\bullet,\bullet}_{z}; d_1, d_2, \delta)$, or rather, 
fix $n$, we study the double complex 
$(\K^{\bullet,n,\bullet}_{z}; d_1, \delta)$. Consider the spectral 
sequence
\[ E^{m,q}_2(\K^{\bullet,n,\bullet}_{z})=
H^m_{d_1} H^q_{\delta}(\K^{\bullet,n,\bullet}_{z}). \]
Since $H^q_{\delta}(\K^{\bullet,n,\bullet}_{z})$ is just 
$H^q(G_z, \LLL^{\bullet,n}_z)$, which is the direct  sum of
subcomplexes of the following form
for all $y\mid \bar z$, $\deg y=-n$:
\begin{equation} \label{complex11} \begin{split}
 0\rightarrow H^q(G_z,& [\BBB_{\mathbf 1},y]) \stackrel{d^1_1}{\rightarrow} 
\cdots \stackrel{d^1_1}{\rightarrow} \\ &
\bigoplus_{\substack{y'\mid y\\ \deg y'=-p}} 
H^q(G_z, [\BBB_{z(y')}, y/y'])
  \cdots\stackrel{d^1_1} {\rightarrow}
H^q(G_z, [\BBB_{z(y)},\mathbf 1]) \rightarrow 0 
\end{split} \end{equation}
where 
\[ [\BBB_{z'},y']:=\langle [a,y']: a\in B_{z'}\rangle_{\T}\cong \BBB_{z'}. \]
Note that for any $y'\mid y\mid \bar z$,
\[ \BBB_{z(y)}=\BBB_{z(y')}\otimes_{\T} \T[G_{z(y/y')}]. \]
One has a commutative diagram
\[ \begin{CD}
H^q(G_z, [\BBB_{z(y')}, y]) @>{-\omega(x, y)d^1_{1,x}}>> 
H^q(G_z, [\BBB_{z(y'/x)},y/x])\\
@V{\theta}VV @V{\theta}VV \\
H^q(G_{z/z(y')}, [\BBB_{\mathbf 1}, y]) @>{res}>>
H^q(G_{z/z(y'x)}, [\BBB_{\mathbf 1},y/x])
\end{CD} \]
where $\theta$ is the isomorphism induced by Shapiro's Lemma. 
Note that $[\BBB_{\mathbf 1}, y]$ is just one copy of $\T$ indexed by
$y$, we write it as $\T[y]$.
Through $\theta$, the complex~\eqref{complex11} is then
quasi-isomorphic to
\begin{equation} \label{complex12} 
 0\rightarrow H^q(G_z, \T[y]) 
\cdots \rightarrow
\bigoplus_{\substack{y'\mid y\\ \deg y'=-p}} 
H^q(G_{z/z(y')}, \T[y/y']) 
\cdots\rightarrow
H^q(G_{z/z(y)}, \T[\mathbf 1]) \rightarrow 0 
\end{equation}
with the differential
\[  \tilde d (c)=-\sum_{x\mid y/y'} \omega(x,y/y') res_{x} c \]
for 
\[ c\in H^q(G_{z/z(y')}, \T[y/y']),\ 
res_{x}\ \text{is the restriction of $c$ in}\
H^q(G_{z/z(xy')}, \T[y/y'x]).  \]  
If replace $q$ in the complex~\eqref{complex12} above by $\ast$, then
we have a complex 
\begin{equation} \label{complex13} 
 0\rightarrow H^{\ast}(G_z, \T[y]) 
\cdots \rightarrow
\bigoplus_{\substack{y'\mid y\\ \deg y'=-p}} 
H^{\ast}(G_{z/z(y')}, \T[y/y']) 
\cdots\rightarrow
H^{\ast}(G_{z/z(y)}, \T[\mathbf 1]) \rightarrow 0 
\end{equation}

\begin{Lemma} \label{Lemma:restriction1}
The complex~\eqref{complex13} is acyclic
except at the first cohomology while the first cohomology is
the direct sum of free graded $\T/m_w\T$-modules $H_{\T,w}$ for 
$y\mid w\mid z$, where $ m_w=\gcd\{|G_{z(x)}|: x\mid w\}$ and
the grading of $H_{T,w}$ is as stated in 
Proposition~\ref{Proposition:trivial1}.
\end{Lemma}
\begin{proof} Since $\T$ is a trivial $G_z$-module, we can apply the
results of Proposition~\ref{Proposition:trivial1} here.
The first cohomology is just
\[ \bigcap_{x\mid y} \ker(H^{\ast}(G_z,\T)\rightarrow 
H^{\ast}(G_{z/z(x)}, \T)), \]
which is nothing but the direct sum of $H_{\T,w}$ for $y\mid w\mid z$
by Proposition~\ref{Proposition:trivial1}.
Apply Proposition~\ref{Proposition:trivial1} 
again, we see the complex~\eqref{complex13} satisfies the conditions of
Lemma 5.2 of Ouyang~\cite{Ouyang2}, Page 16. Following that lemma, we 
know  other cohomology groups vanish for the complex~\eqref{complex13}.
\end{proof}
Write $H_{\T,w}^q$ the $q$-th component of $H_{\T,w}$, we thus have
\begin{Proposition}
For any fixed $n$, 
the $E^{m,q}_2$ term $H^m_{d_1} H^q_{\delta} 
(\K^{\bullet,n,\bullet}_{z})$ of the double complex
$(\K^{\bullet,n,\bullet}_{z}; d_1,\delta)$,
is then the direct sum of free 
$\T/m_{w}\T$-modules  $H^q_{\T,w}[y]$  where
\[ \deg y=-n,\qquad  y\mid \bar w\mid z \]
and the $\T/m_{w}\T$-rank of $H^q_{\T,w}[y]$ is $\binom{\deg \bar w-1}{i-1}$
if $q=2\deg w-\deg \bar w+i$.
\end{Proposition}

\subsection{The case $p(x;1)=0$ for every $x\mid z$} 
\label{subsection:quotient}
In this subsection, we suppose that $p(x;1)=0$ for every $x\mid z$.
In this case, we can give a complete description of the
$G_z$-cohomology of $\NNN_z$. 
Consider the $\T$-submodule 
$\mathbf S$ of $\K_{z}$ generated by 
\[ \{[a, y, w]: a\in B_{z/z(y)},\ y\mid z,\ 
\bar w\mid z, a\notin B_{\mathbf 1}\ \text{if}\ y\mid w \}.
\]
Under the assumption $p(x;1)=0$, 
one easily sees that $d_1 \mathbf S, d_2 \mathbf S, \delta \mathbf S\subseteq 
\mathbf S$, thus $\mathbf S$ is really a subcomplex of $\K_{z}$ with 
related double and triple complex structures. We let $\Q_{z}
=\K_{z}/\mathbf S$, thus $\Q_{z}$ is a free $\T$-module generated by
\[ \{[\mathbf 1, y, w]: y\mid \bar w\mid z \}. \] 
Note that the induced differential $d_1=0$ in $\Q_{z}$. We write 
the quotient map as $\rho$.

\begin{Proposition} \label{Proposition:quasi1}
The quotient map $\rho$ is a quasi-isomorphism.
\end{Proposition}

\begin{proof} Consider the  triple complex $(\Kc_{z}; d_1, d_2,\delta)$
and the related triple complex $(\Q^{\bullet,\bullet,\bullet}_{z}; d_1,
d_2, \delta)$. Fix $d_2$-degree $n$, we consider the double complex
$(\K^{\bullet,n,\bullet}_{z}; d_1,\delta)$ and its quotient by $\rho$.
Then $\rho$ induces a map  
\[ \rho_2: H^m_{d_1}(H^q_{\delta}(\K^{\bullet,n,\bullet}_{z}))
\longrightarrow
H^m_{d_1}(H^q_{\delta}(\Q^{\bullet,n,\bullet}_{z})). \]
We claim that $\rho_2$ is an isomorphism.

Assuming the claim, then
$H^{m+q}_{total}(\K^{\bullet, n,\bullet}_z, d_1+\delta)$ is isomorphic to 
$H^{m+q}_{total} (\Q^{\bullet, n,\bullet}_z, d_1+\delta)$. Thus for the double 
complex $(\K^{\bullet,\bullet}_{z}; d_2, d_1+\delta)$
and its quotient   $(\Q^{\bullet,\bullet}_{z}; d_2, d_1+\delta)$, 
the $E_2^{n, m+q}$-term
$H^{m+q}_{d_2}(H^n_{d_1+\delta}(\K^{\bullet,\bullet}_{z}))$ is 
isomorphic to
$H^{m+q}_{d_2}(H^n_{d_1+\delta}(\Q^{\bullet,\bullet}_{z}))$. 
$\rho$ hence is a quasi-isomorphism. Noe that here we use the 
following fact about spectral sequences: a complex homomorphism 
is a quasi-isomorphism if in the corresponding 
weakly convergent spectral sequences,
the $E_r$-terms are isomorphic for some positive integer $r$.

Now we show the isomorphism of 
$\rho_2$. Consider the complex $(L^{\bullet}_{y}, \delta)$ 
generated by $\{ [\mathbf 1, y, w]: \bar w\mid z \}$. This 
complex is exactly $\Hom(P_{z\bullet}, [\BBB_{\mathbf 1},y])$.
Let $L^{\prime\bullet}_{y}$ and $L^{\prime\prime\bullet}_{y}$ be 
the subcomplexes
generated by $\{[\mathbf 1,y, w]: y\mid w\}$ and by 
$\{[\mathbf 1, y, w]: y \nmid w\}$ respectively. 
Thus $L^{\bullet}_{y}$ is the direct sum of
$L^{\prime\bullet}_{y}$ and $L^{\prime\prime\bullet}_{y}$. 
Correspondingly, $H^{\ast}(G_z,[\mathbf 1, y])$ is the direct sum of
$H^{\ast}(L^{\prime\bullet}_{y}, \delta)$ and
$H^{\ast}(L^{\prime\prime\bullet}_{y}, \delta)$. Now
the kernel of $d^1_1$ at $H^{q}(G_z,[\mathbf 1, y])$
in the complex~\eqref{complex11}, or equivalently,
in the complex~\eqref{complex12}, is just 
$H^{q}(L^{\prime\bullet}_{y}, \delta)$. We see that 
$\Q^{\bullet,n,\bullet}_{z}$ is actually the direct sum of 
$L^{\prime\bullet}_{y}$(Note that $d_1=0 $ in $\Q_z$). This proves 
the isomorphism of $\rho_2$.
\end{proof}

\begin{theo} \label{theo:b}
If for every $x\mid z$, $p(x;1)=0$. Then 
$H^{\ast}(G_z,\NNN_z)$, the $G_z$-cohomology
of the universal norm distribution $\NNN_z$ is the direct sum of 
$H_{\T,w}[y]$ where $H_{\T,w}$ is as stated in 
Proposition~\ref{Proposition:trivial1} and
\[ y\mid \bar w\mid z. \]
Any element 
$c[y]\in H^{\ast}(G_z,\NNN_z)$ for $c\in H^q_{\T,w}$ is of degree
$q-\deg y$.
\end{theo}

\begin{rem} Let $\NNN_z$=$U_r$, the universal ordinary distribution of
level $r$, if $r$ is odd, then $G_{p^i}$ is cyclic for every $p^i\| r$.
We also see that $p(x;1)=1-1=0$, hence the above theorem gives a
complete description of $H^{\ast}(G_r,U_r)$ and generalizes 
Theorem A
in Ouyang~\cite{Ouyang2}, where we need the condition
$r$ is squarefree.
\end{rem}

\subsection{The $G_z$-cohomology of $\NNN_z/M\NNN_z$} We suppose now that
$M$ is a common divisor of $|G_{z(x)}|$ and $p(x;1)$ for every $x\mid z$.
Let $\mathbf S_z$ be the same as in \S~\ref{subsection:quotient}. Then
$\mathbf S_z/M\mathbf S_z$ is a submodule of $\K_{z}/M\K_z$ generated by 
\[ \{[a, y, w]: a\in B_{z/z(y)},\ y\mid z,\ 
\bar w\mid z, a\notin B_{\mathbf 1}\ \text{if}\ y\mid w \}.
\]
One easily sees that $\mathbf S_z/M\mathbf S_z$ is a subcomplex of 
$\K_{z}/M\K_z$ with respect to the multi-complex structure of $\K_z/M\K_z$. 
We let $\Q_{z}/M\Q_z$ be the quotient of $\K_z/M\K_z$ to $\mathbf S_z/
M\mathbf S_z$, thus $\Q_{z}/M\Q_z$ is a free $\T/M\T$-module generated by
\[ \{[\mathbf 1, y, w]: y\mid \bar w\mid z \}. \] 
Note that the induced differentials $d_1=d_2=d=\delta=0$ in $\Q_{z}/M\Q_z$. 
Write 
the quotient map from $\K_z/M\K_z$ to $\Q_z/M\Q_z$ as $\rho_M$.
\begin{Proposition} \label{Proposition:quasi2}
The homomorphism $\rho_M$ is a quasi-isomorphism.
\end{Proposition}
\begin{proof}
Similar to the proof of Proposition~\ref{Proposition:quasi1}
\end{proof}

\begin{theo} \label{theo:a}
Let $M\in \co$ be a common divisor of $|G_{z(x)}|$ and
$p(x;1)$ for all $x\mid z$. Then
the cohomology group $H^{\ast}(G_z,\NNN_z/M\NNN_z)$ is a direct sum 
of rank one graded $\T/M\T$-modules $\langle c(y,w)\rangle$ where
\[ y\mid \bar w\mid z, \deg c(y,w)=\deg w-\deg y. \]
\end{theo}
\begin{proof}
By the quasi-isomorphism of $\rho_M$ in Proposition~\ref{Proposition:quasi2},
the cohomology group $H^{\ast}(G_z,\NNN_z/M\NNN_z)$ is then just the 
total cohomology group of the complex $\Q_z/M\Q_z$. However, all induced
differentials in $\Q_z/M\Q_z$ are $0$, thus its cohomology is itself. Let
$c(y,w)$ be the element in $H^{\ast}(G_z,\NNN_z/M\NNN_z)$ represented by
the cocycle $[\mathbf 1,y,w]$ in $\Q_z/M\Q_z$, we hence get the proof
of the above theorem.
\end{proof}

\begin{rem} \label{rem:theoa}
With the automorphism $\epsilon$ in \S~\ref{subsection:another}, 
we easily see that  
\[ \rho_M: (\K^{\bullet,\bullet}_z/M\K^{\bullet,\bullet}_z;
\tilde d,\tilde \delta)\rightarrow 
(\Q^{\bullet,\bullet}_z/M\Q^{\bullet,\bullet}_z;0,0) \]
is a quasi-isomorphism, thus Theorem~\ref{theo:a} can be stated 
in the form of the double complex $(\K^{\bullet,\bullet}_z;
\tilde d,\tilde \delta)$.

We call the basis $\{c(y,w): y\mid \bar w\mid z\}$ given in 
Theorem~\ref{theo:a}, 
the {\em canonical basis} for   $H^{\ast}(G_z, 
\NNN_z/M\NNN_z)$. In particular, we write $c(y,y)$ as $c_y$.
By the above theorem, we see that for every $z\in Z$,
\[ H^0(G_z, \NNN_z/M\NNN_z)=\langle c_y:\ y\mid z\rangle_{\T/M\T} \]
is the union of all
$H^0(G_{z'}, \NNN_{z'}/M\NNN_{z'})$ with $z'\mid_s z$ and $z'$ finite.
We'll use this fact in Ouyang~\cite{Ouyang4}
for the double complex $(\K^{\bullet,\bullet}_z;
\tilde d,\tilde \delta)$.
\end{rem}

\begin{rem}
One can expect parallel result to Theorem B in Ouyang~\cite{Ouyang2}
holds here too. The answer is yes. However, we feel more
appropriate to state it in Ouyang~\cite{Ouyang4}, as a natural 
consequence of the universal Kolyvagin recursion, just like
the proof of the above Theorem B in Anderson and Ouyang~\cite{AO}.
\end{rem}

\end{document}